\documentclass[a4paper]{article} 
\usepackage{graphicx} 
\usepackage{geometry} 
\usepackage[nottoc]{tocbibind} 
\usepackage{enumitem} 

\usepackage{color} 
\definecolor{newcolor}{rgb}{0.0, 0.6, 0.0}
 
\usepackage{color} 
\definecolor{newcolor2}{rgb}{0.0, 0.0, 0.8}
 
\usepackage[utf8]{inputenc}
\usepackage{mathrsfs}
\usepackage{amssymb}
\usepackage{amsmath} 
\usepackage{amsthm} 
\usepackage{faktor}
\usepackage{dsfont}
\usepackage{notoccite}
\usepackage{tikz}
\usepackage{caption}
\usepackage{subcaption}
\newcommand{\R}{\mathbb{R}}

\newcommand{\E}{\mathbb{E}}

\newcommand{\1}{\mathds{1}}
\newcommand{\s}{\hspace{1pt}}
\newcommand{\Var}{\text{Var}}

\newcommand{\Ent}{\text{Ent}}
\newtheorem{theorem}{Theorem}[section]

\newtheorem{example}{Example}[section]

\title{A Fuzzy Decomposition Method to Establish Functional Inequalities}
\author{Suei-Wen Chen\thanks{This work is part of an undergraduate thesis done under the supervision of Prof. Pierre Youssef.} \\
\normalsize{Department of Science, New York University Abu Dhabi, UAE}}

\date{}

\begin{document}

\maketitle 

\begin{abstract}
   This paper presents a method to establish functional inequalities via fuzzy decomposition on the state space, which generalizes earlier results dealing with exact partitions of the state space. Given a reversible Markov chain on a finite state space, we define its projection chain and restriction chains from classes of a fuzzy partition on the state space. The Poincaré, log-Sobolev and modified log-Sobolev inequalities associated with the original chain can be estimated from those of its projection chain and restriction chains which tend to have simpler structures and hence easier to work with. An application of this generalization is presented in which the fuzzy decomposition method applies but its exact partitioning counterpart does not.
\end{abstract}

\tableofcontents

\section{Introduction}
Finite-state Markov chains are used in many fields to describe and predict the dynamics of sequential processes. One powerful technique to study the behavior of Markov chains is via functional inequalities (for instance, see \cite{princeton-high-dim-prob, mlsi, lsi}). Applications of functional inequalities includes establishing concentration of measure as well as providing non-asymptotic upper bounds on the mixing time of the underlying Markov chains, which is of importance for applications such as Markov Chain Monte Carlo algorithms (for example, see \cite{markov-chains-and-mixing-times}). Although functional inequalities are powerful in theory, it is often difficult to establish sharp bounds for these inequalities in practice. This paper provides a generalization of a method for establishing functional inequalities via Markov chain decomposition.

\subsection{Reversible Markov Chains and Functional Inequalities}
This paper concerns reversible Markov chains on finite state spaces. Consider a finite probability space $(\Omega, \mathcal{P}(\Omega), \pi)$ and a Markov generator $Q$ with respect to $\pi$, which is an $|\Omega|\times|\Omega|$ matrix with nonnegative off-diagonal entries such that each row sums to zero. We require that $Q$ is reversible under the probability distribution $\pi$, namely $\pi(x)Q(x,y)=\pi(y)Q(y,x)$ for all $x,y\in\Omega$. It is well-known that if $Q$ is irreducible, then the continuous-time chain defined by $p_t:=e^{tQ}$ converges to the stationary distribution, namely $p_t(x,y)\to \pi(y)$ as $t\to\infty$ for all $x,y\in \Omega$. The mixing time defined by
\[t_{\text{mix}}(\epsilon) := \inf \{t\geq 0: \max_{x\in\Omega} \|p_t(x,\cdot)-\pi\|_{\text{TV}}\leq \epsilon \}\]
measures the rate at which this convergence occurs, where $\|\mu-\nu\|_{\text{TV}}=\max_{A\subseteq\Omega} |\mu(A)-\nu(A)|$ is the total variation distance and  $\epsilon>0$ is the threshold.

We denote by $\E_\pi [\cdot]$, $\Var_\pi (\cdot)$ and $\Ent_\pi (\cdot)$ the standard expectation, variance and entropy over the probability space $(\Omega, \mathcal{P}(\Omega), \pi)$, namely
\begin{align*}
\E_\pi [f] & := \sum_{x\in\Omega}\pi(x)f(x) ; \\
\Var_\pi (f) & :=\E_\pi[f^2]-(\E_\pi [f])^2 ; \\
\Ent_\pi (f) & :=\E_\pi [f\log f] -\E_\pi [f]\log \E_\pi [f].
\end{align*}
Variance and entropy capture the global variation of $f$ under $\pi$, while the local variation associated with $Q$ can be measured by the Dirichlet form
\begin{align*}
    \mathcal{E}(f,g) :=\langle -Qf, g\rangle_\pi = -\sum_{x\in\Omega}(Qf)(x)g(x)\pi(x)
\end{align*}
for $f,g:\Omega\to\R$. Establishing bounds on global variations using average local variations provides powerful control over the mixing behavior of the chain. Specifically, the Poincaré constant $\lambda(Q)$, the modified log-Sobolev constant $\alpha(Q)$ and the log-Sobolev constant $\rho(Q)$ are defined respectively as the supremum of $\lambda, \alpha, \rho \geq 0$ such that 
\begin{align*}
    \mathcal{E} (f,f) \geq \lambda \Var_\pi(f);  \\
    \mathcal{E}(f,\log f) \geq \alpha \Ent_\pi(f); \\
    \mathcal{E}(\sqrt{f},\sqrt{f}) \geq \rho \Ent_\pi(f), 
\end{align*}
for all $f:\Omega\to (0,\infty)$. The three inequalities are called the Poincaré inequality, the modified log-Sobolev inequality  (MLSI) and the log-Sobolev inequality (LSI) respectively. It can be shown that the mixing time is of order $O(\frac{1}{\lambda(Q)}\log \pi_{\min}^{-1})$, $O(\frac{1}{\alpha(Q)}\log\log \pi_{\min}^{-1})$ and $O(\frac{1}{\rho(Q)}\log \log \pi_{\min}^{-1})$, where $\pi_{\min}:=\min_{x\in\Omega} \pi(x)$ and the three constants satisfy $2\lambda(Q) \geq \alpha(Q)\geq 4\rho(Q)$ (see \cite{lsi,mlsi}). 

\subsection{Decomposition Methods and Fuzzy Decomposition}
Powerful as these functional inequalities are, proving sharp lower bounds for these constants remains difficult in practice. In fact, proving the optimal LSI and MLSI for even a two-state chain is highly nontrivial \cite{conc-of-meas-and-lsi}. It is thus desirable to find a reduction method for establishing such functional inequalities. A natural way to proceed is to partition the state space $\Omega$ into disjoint subsets. This naturally induces several restriction chains, one on each block, as well as a projection chain between blocks. The restriction and projection chains could be simpler in structure and easier to establish their corresponding functional inequalities. The task then is to recover the functional-inequality constant for the original chain from those of the restriction chains and the projection chain. A review of such decomposition methods can be found in \cite{decomposable-mc}.

The main contribution of this paper is a generalization of a partition method developed by Lu and Yau in the context of interacting particle systems \cite{partition-original}. A succinct discussion on this method is given by Hermon and Salez in \cite{partition}, which builds on the usual notion of partition. Specifically, let $I$ be a finite set whose elements are thought of as classes to partition $\Omega$ into. Traditionally, a partition on $\Omega$ is given by pairwise disjoint nonempty subsets $\Omega_i\subseteq\Omega$ such that $\Omega=\cup_{i\in I}\Omega_i$. We refer to this as exact partition in this paper. In some cases, however, it may be natural to relax the condition of exact membership and interpret each state as belonging to a number of classes with different degrees of membership. This motivates the definition of \emph{fuzzy decomposition}, which we describe below in the language of fuzzy sets and fuzzy relations \cite{fuzzy-set}.

Given a set $X$, a fuzzy set on $X$ is defined as a map $\mu:X\to [0,1]$. A subset $A\subseteq X$ can be thought of as the collection of elements in $X$ satisfying some property $P$. A fuzzy set $\mu$ on $X$ is meant to model the degree to which each $x\in X$ satisfies some property $P$. The closer $\mu(x)$ is to one (resp. zero), the more (resp. less) we consider $x$ to satisfy $P$. A fuzzy relation between sets $X$ and $Y$ is a fuzzy set on $X\times Y$. The underlying mapping $X\times Y \to [0,1]$ quantifies the degree to which this relation is satisfied.

We define a fuzzy decomposition (or a fuzzy partition) of $\Omega$ with classes $I$ as a fuzzy relation $a:\Omega\times I \to [0,1]$ on $\Omega\times I$ such that (1) $\sum_{i\in I} a(x,i)=1$ for all $x\in\Omega$ and (2) for all $i\in I$ we have $a(x,i)>0$ for some $x\in \Omega$. We write $a_i(x):= a(x,i)$ and interpret this as the degree to which $x$ belongs to class $i$. The first condition means that for each state $x$, the degrees of membership across all classes sum up to unity, while the second condition rules out empty classes. For each $i\in I$ we refer to the set $\Lambda_i := \{x\in\Omega: a_i(x)>0\}$ as class $i$, consisting of states with positive membership value. The traditional notion of partition is recovered by restricting the range of $a$ to $\{0,1\}$ and set $\Omega_i$ to $\Lambda_i$. The main result of the paper generalizes the treatment in \cite{partition} to the case of fuzzy decompositions.

\subsection{Main Result}

Suppose we are given a reversible chain $(Q,\pi)$ on $\Omega$ and a fuzzy decomposition $a$ on $\Omega$ with classes $I$, where
\[\Lambda_i:=\{x\in\Omega \mid a_i(x)>0\} \]
is referred to as class $i$. The stationary measure $\pi$ yields a probability measure $\hat{\pi}$ on $I$ defined by
\[\hat{\pi}(i) := \sum_{x\in \Lambda_i}a_i(x)\pi(x), \s\s\s i\in I\]
as well as a probability measure on each class $\Lambda_i$ defined by
\[\pi_i(x) := \frac{a_i(x)\pi(x)}{\hat{\pi}(i)}, \s\s\s x \in \Lambda_i.\]
Define the projection chain $\hat{Q}$ of $Q$ as follows: for $i, j\in I$, $i\neq j$ let
\[\hat{Q}(i,j) := \frac{1}{\hat{\pi}(i)}\sum_{x\in\Lambda_i}\sum_{y\in\Lambda_j\setminus\{x\} }a_i(x)a_j(y)\pi(x)Q(x,y),\]
and set $\hat{Q}(i,i)$ such that each row of $\hat{Q}$ sums to $0$. For each $i\in I$, the restriction chain $Q_i$ of $Q$ on $\Lambda_i$ is defined as follows: $Q_i(x,y):=a_i(y)Q(x,y)$ if $x\neq y$ and set $Q_i(x,x)$ such that each row sums to $0$. 

It is clear from their definition that $\hat{Q}$ is reversible under $\hat{\pi}$ and that $Q_i$ is reversible under $\pi_i$ for all $i \in I$. The interpretations of $\hat{Q}$ and $Q_i$ are natural: $\hat{Q}(i,j)$ represents the transition probability from class $i$ to class $j$ as a normalized sum of transition probabilities from states in class $i$ to states in class $j$, weighted by proportions of class membership. The transition probability $Q_i(x,y)$ is simply the original transition probability $Q(x,y)$ discounted by $a_i(y)$, namely how strongly $y$ belongs to class $i$.

Now suppose that for $(i,j)\in I^2$ with $\hat{Q}(i,j)>0$ we are given couplings $\kappa_{ij}$ of $\pi_i$ and $\pi_j$, meaning that $\kappa_{ij}$ is a probability distribution on $\Lambda_i\times\Lambda_j$ such that for all $x\in\Lambda_i$ and $y\in\Lambda_j$,
\[\sum_{x\in\Lambda_i} \kappa_{ij}(x,y)=\pi_j(y) \text{ and } \sum_{j\in\Lambda_j} \kappa_{ij}(x,y)=\pi_i(x). \] We define the quality of these couplings to be
\[ \chi := \min_{\substack{(x,y,i,j) \text{ s.t.}\\ \kappa_{ij}(x,y)>0,\\ x\neq y}} \dfrac{a_i(x)a_j(y)\pi(x)Q(x,y)}{\hat{\pi}(i)\hat{Q}(i,j)\kappa_{ij}(x,y)} .\]
This quantity allows us to estimate the Poincare, log-Sobolev and modified log-Sobolev constants of $Q$ from those of $Q_i$ and $\hat{Q}$ in the following manner.

\begin{theorem}\label{fuzzy-partition} \normalfont
Let  $(\Omega, \mathcal{P}(\Omega), \pi)$ be a finite probability space. Suppose we are given a fuzzy partition on $\Omega$ with classes $I$ which admits couplings $\kappa_{ij}$ between $\pi_i$ and $\pi_j$ for all $(i,j)\in I^2$ such that $\hat{Q}(i,j)>0$, with quality $\chi$. Then
\begin{align*}
    \lambda(Q) &\geq \min \{ \chi\lambda(\hat{Q}), \min_{i\in I} \lambda(Q_i)\}; \\
    \alpha(Q) &\geq \min \{ \chi\alpha(\hat{Q}), \min_{i\in I} \alpha(Q_i)\}; \\
    \rho(Q) &\geq \min \{ \chi\rho(\hat{Q}), \min_{i\in I} \rho(Q_i)\}.
\end{align*}
\end{theorem}

Let us note that Theorem \ref{fuzzy-partition} is a strict generalization of its non-fuzzy version in \cite{partition-original, partition}. This is because in order for such decomposition methods to give nontrivial bounds, the projection chain and all the restriction chains must be irreducible, and at the same time the existence of couplings with nonzero quality is needed. For some Markov chains, it is not always possible to achieve all these requirements simultaneously. For instance, the non-fuzzy version of Theorem \ref{fuzzy-partition} cannot be applied to Example \ref{graph-example} in general due to the asymmetry of the state space.

\section{Related Works and Examples}
A number of decomposition methods have been developed for  Markov Chain analysis. In the method proposed by Madras and Randall \cite{mc-decomp-sampling}, decomposition takes the form of a set cover $\Omega=\bigcup_{i\in I} A_i$  where $A_i\subseteq\Omega$ are not necessarily pairwise disjoint. A different version of projection chain and restriction chains is considered. In particular, the projection chain depends only on the probabilities of $A_i$ and the maximum overlap, defined as the maximum cardinality of $\{i\in I : x\in A_i\}$ across $x\in\Omega$. The method by Jerrum et al. \cite{decomposable-mc} works with exact partition $\Omega=\bigcup_{i\in I} \Omega_i$ of the state space. While their version of projection chain and restriction chain coincide with the one considered in this paper, they do not appeal to couplings of the distributions $\pi_i$. Instead, the quantity
\[\gamma:= \max_{i\in I} \max_{x\in\Omega_i} P(x, \Omega\setminus\Omega_i)\]
describing the maximum probability of exiting the current block  is considered, where $P$ is the transition probability matrix of the chain. However, unlike Theorem \ref{fuzzy-partition}, these two settings are either insensitive to or unable to capture the degree of membership that may be desirable in some applications. 

The following is an example of how Theorem \ref{fuzzy-partition} may apply.

\begin{example} \normalfont \label{graph-example}
Let $G$ be a finite set on which an undirected connected graph is defined. We consider a graph consisting of two copies of $G$ but with some points glued together. More precisely, let $H$ be a subset of $G$ such that if $x\in H$ then $N_G(x)\cap H=\varnothing$, where $N_G(x)$ is the set of neighbors of $x$ in $G$. Define an equivalence relation $\equiv$ on $G\times\{1,2\}$ by $(x,i)\equiv (y,j)$ if and only if $(x,i)= (y,j)$ or $x=y\in H$. Let $\mathcal{G}:=\faktor{G\times\{1,2\}}{\equiv}$. We identify the elements of $H$ with their equivalent classes; namely, since $(h,1)\equiv(h,2)$ for $h\in H$ we abuse the notation and write $H=\{(h,i): h\in H, i\in \{1,2\}\}$. We also write $G_i=\{(x,i)\in\mathcal{G}\mid x\in G\setminus H\}$. Thus $\mathcal{G}=G_1\cup G_2\cup H$. The edges of $\mathcal{G}$ are defined by the edges of the two copies of $G$ while connecting $(x,1)$ and $(x,2)$ for all $x\in G\setminus H$. Write $v\sim w$ if $v$ and $w$ are neighbors in $\mathcal{G}$. See Figure \ref{graph-example}
for an illustration of this construction for $G=\{a,b,c,d,e\}$ with edges $a\sim b$, $b\sim c$, $c\sim d$, $d\sim e$, $e\sim b$, $e\sim a$ with $H=\{a,c\}$.
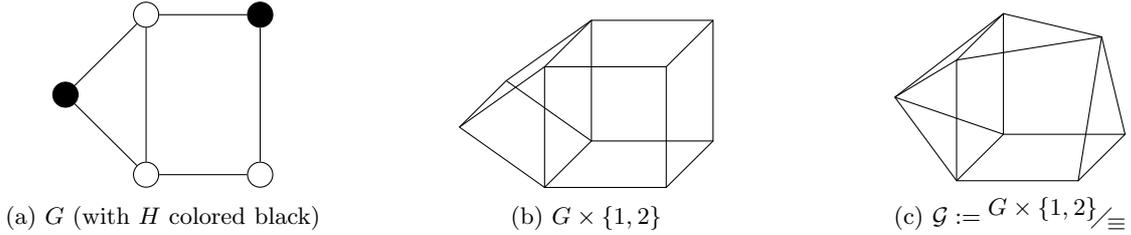
\begin{figure}
     \centering
     \begin{subfigure}[b]{0.3\textwidth}
         \centering
\begin{tikzpicture}[main/.style = {draw, circle}, blacknode/.style={circle, draw=black, fill=black}, node distance=15mm] 
            \node[blacknode] (h){}; 
            \node[main] (1) [above right of=h]{};
            \node[main] (2) [below right of=h]{};
            \node[main] (3) [right of=2]{};
            \node[blacknode] (4) [right of=1]{};
            \draw (h)--(1);
            \draw (h)--(2);
            \draw (1)--(2);
            \draw (1)--(4);
            \draw (2)--(3);
            \draw (3)--(4);
        \end{tikzpicture}
         \caption{$G$ (with $H$ colored black)}
         \label{fig:G}         
     \end{subfigure}
     \hfill
     \begin{subfigure}[b]{0.3\textwidth}
         \centering
\begin{tikzpicture}[scale=0.8]
            \draw (0,0,0) -- (2,0,0) -- (2,2,0) -- (0,2,0) -- cycle;
            \draw (0,0,2) -- (2,0,2) -- (2,2,2) -- (0,2,2) -- cycle;
            \draw (0,0,0) -- (0,0,2);
            \draw (2,0,0) -- (2,0,2);
            \draw (0,2,0) -- (0,2,2);
            \draw (2,2,0) -- (2,2,2);
            \draw (0,0,0) -- (-1.4,1,0);
            \draw (0,2,0) -- (-1.4,1,0);
            \draw (0,0,2) -- (-1.4,1,2);
            \draw (0,2,2) -- (-1.4,1,2);
            \draw (-1.4,1,0) -- (-1.4,1,2);
    \end{tikzpicture}          
         \caption{$G\times\{1,2\}$}
         \label{fig:two-G}
     \end{subfigure}
     \hfill
     \begin{subfigure}[b]{0.3\textwidth}
         \centering
\begin{tikzpicture}[scale=0.8]
            \draw (0,0,0) -- (2,0,0) -- (2,2,1) -- (0,2,0) -- cycle;
            \draw (0,0,2) -- (2,0,2) -- (2,2,1) -- (0,2,2) -- cycle;
            \draw (0,0,0) -- (0,0,2);
            \draw (2,0,0) -- (2,0,2);
            \draw (0,2,0) -- (0,2,2);
            \draw (0,0,0) -- (-1.4,1,1);
            \draw (0,2,0) -- (-1.4,1,1);
            \draw (0,0,2) -- (-1.4,1,1);
            \draw (0,2,2) -- (-1.4,1,1);
    \end{tikzpicture}          
         \caption{$\mathcal{G}:=\faktor{G\times\{1,2\}}{\equiv}$}
         \label{fig:five over x}
     \end{subfigure}
\caption{Illustration of the construction of $\mathcal{G}$}
\label{graph-example}
\end{figure}

Now we consider the simple random walk on $\mathcal{G}$. Let $d_G(x)$ and $d_\mathcal{G}(x)$ denote the degree of $x$ in $G$ and $\mathcal{G}$ respectively, and let
\[D_G:=\sum_{x\in G}d_G(x), \hspace{5pt} D_\mathcal{G} := \sum_{x\in\mathcal{G}} d_\mathcal{G}(x).\]
Then the chain in question is $Q(x,y)=1/d_\mathcal{G}(x)$ if $x\sim y$ and $\pi(x)=d_\mathcal{G}(x) / D_\mathcal{G}$, where
\begin{equation*}
    d_\mathcal{G}(x)= \begin{cases}
    2d_G(x) &\text{if } x\in H, \\
    d_G(x)+1 &\text{otherwise.}
\end{cases}
\end{equation*}
Thus $D_\mathcal{G}=2D_G+2\#(G\setminus H)$ where $\#$ denotes the cardinality of a set. Assign weight vectors of $x\in\mathcal{G}$ to be $a(x)=(1/2,1/2)$ if $x\in H$, $a(x)=(1, 0)$ if $x\in G_1$ and $a(x)=(0,1)$ if $x\in G_2$. Thus $\Lambda_i=G_i\cup H$ and $\hat{\pi}(i)=1/2$  by symmetry ($i=1,2$). Note that $\pi_1(x)=\pi_2(x)=\pi(x)$ for $x\in H$ and $\pi_1(x)=\pi_2(y) = 2\pi(x)$ if $x\in G_1$ and $y\in G_2$ correspond to the same element in $G$. Hence the probability measure $\kappa:=\kappa_{12}$ on $\Lambda_1\times \Lambda_2$  defined by
\begin{equation*}
    \kappa(x,y)=\begin{cases}
        \pi(x) &\text{if } x=y\in H \\
        2\pi(x) &\text{if } x=(v,1) \text{ and } y=(v,2) \text{ for some } v\in G\setminus H \\
        0 & \text{otherwise}
    \end{cases}
\end{equation*}
is a coupling of $\pi_1$ and $\pi_2$. Again, by symmetry we set $\kappa_{21}(y,x)=\kappa(x,y)$. Note that $\pi(x)Q(x,y)=\1_{x\sim y}/D_\mathcal{G}$, so the projection chain is
\begin{align*}
    \hat{Q}(2,1) & = \hat{Q}(1,2)= \frac{2}{D_\mathcal{G}} \left( \sum_{x\in H} d_G(x) + \#(G\setminus H) \right).
\end{align*}
Recall that the quality of the couplings are
\begin{align*}
    \chi = \min \frac{a_1(x)a_2(y)\pi(x)Q(x,y)}{\hat{\pi}(1)\hat{Q}(1,2)\kappa(x,y)}
\end{align*}
where the minimum is taken for $x\neq y$ with $\kappa(x,y)>0$ and therefore $x=(v,1)$, $y=(v,2)$ for some $v\in G\setminus H$. It follows that 
\begin{align*}
    \chi = \min_{x\in G_1\cup G_2} \frac{1/D_\mathcal{G}}{\frac{1}{2} \hat{Q}(1,2) \cdot 2d_\mathcal{G}(x)/D_\mathcal{G}} =  \frac{1}{\hat{Q}(1,2)} \frac{1}{\max_{v\in G\setminus H}d_G(v)+1}.
\end{align*}
The restriction chains $Q_i(x,y)$ on $G$ coincide with $Q(x,y)$ when $y\in G\setminus H$ and $Q_i(x,y)=\frac{1}{2}Q(x,y)$ if $y\in H$. We conclude that
\begin{align*}
    \lambda(Q) \geq \min \left\{\frac{2}{\max_{v\in G\setminus H}d_G(v)+1}, \lambda(Q_i) \right\}
\end{align*}
and similarly for $\alpha(Q)$ and $\rho(Q)$. A lower bound on $\lambda(Q)$ can be established if one can solve for $\lambda(Q_i)$, for example, by recursively applying decomposition techniques or spectral methods, depending on the structure of the original graph.

\end{example}

\section{Proof}
A proof of Theorem \ref{fuzzy-partition} is presented in this section, which directly generalizes the treatment in \cite{partition}. First, we note that by reversibility, the Dirichlet forms associated with the three constants can be written as
\begin{equation*}  \label{local-expression}
    \mathcal{L}_\pi f:=\frac{1}{2}\sum_{(x,y)\in\Omega^2} \pi(x)Q(x,y)\Psi(f(x),f(y)) \text{, where}
\end{equation*}  
\begin{equation*} \label{convex-func}
    \Psi(u,v) = \begin{cases}
    (u-v)^2 &\text{for Poincare};\\
    (u-v)(\log u-\log v) &\text{for modified log-Sobolev};\\
    (\sqrt{u}-\sqrt{v})^2 &\text{for log-Sobolev}.
    \end{cases}
\end{equation*}
For notational convenience we also write
\begin{equation*}
    \mathcal{R}_\pi f := \begin{cases}
    \Var_\pi f &\text{for Poincare}\\
    \Ent_\pi f &\text{for modified log-Sobolev and log-Sobolev.}
    \end{cases}
\end{equation*}

The idea is to decompose the variance, entropy and Dirichlet forms according to the partition. To this end, for $i\in I$ we define 
\[\hat{f}(i):=\E_{\pi_i} \left[ \left.f\right\vert_{\Lambda_i}\right] = \sum_{x\in\Lambda_i}f(x)\pi_i(x).\]
Thus
\begin{align*}
    & \Var f = \frac{1}{2}\sum_{x,y\in\Omega}\pi(x)\pi(y)(f(x)-f(y))^2 \\
    & =  \frac{1}{2}\sum_{i,j\in I}\sum_{x,y\in\Omega}a_i(x)a_j(y)\pi(x)\pi(y)[(f(x)-\hat{f}(i)) + (\hat{f}(i)-\hat{f}(j)) + (\hat{f}(j)-f(y)) ]^2.
\end{align*}
When expanding the square, the cross terms involving $(f(x)-\hat{f}(i))$ or $(\hat{f}(j)-f(y))$ are $0$ by the definition of $\hat{f}$. Thus, $\Var_\pi f$ is equal to
\begin{align*}
     & \frac{1}{2}\sum_{i,j\in I}\sum_{x,y\in\Omega}a_i(x)a_j(y)\pi(x)\pi(y)[(f(x)-\hat{f}(i))^2 + (\hat{f}(i)-\hat{f}(j))^2 + (\hat{f}(j)-f(y))^2 ] \\
    & = \sum_{i\in I}\sum_{x\in\Omega}a_i(x)\pi(x)(f(x)-\hat{f}(i))^2 + \frac{1}{2}\sum_{i,j\in I}\hat{\pi}(i)\hat{\pi}(j) (\hat{f}(i)-\hat{f}(j))^2 \\
    & = \sum_{i\in I}\hat{\pi}(i)\Var_{\pi_i} f + \Var_{\hat{\pi}} \hat{f}.
\end{align*}
Similarly, the entropy of $f:\Omega\to (0,\infty)$ can be written as
\begin{align*}
    \Ent_{\pi} f & = \E_{\pi} [f(\log f- \log \E_\pi f)]  = \sum_{x\in\Omega}\sum_{i\in I} a_i(x)\pi(x)f(x)\left[\log f(x) - \log \E_\pi f \right] \\
    & = \sum_{i\in I} \hat{\pi}(i) \sum_{x\in\Lambda_i} \pi_i(x)f(x)\left[\log f(x) - \log \E_\pi f \right] \\
    & = \sum_{i\in I} \hat{\pi}(i)\Ent_{\pi_i} f + \sum_{i\in I}  \hat{\pi}(i) \sum_{x\in\Lambda_i} \pi_i(x)f(x)\left[\log\E_{\pi_i}f - \log \E_\pi f \right] \\
    & = \sum_{i\in I} \hat{\pi}(i)\Ent_{\pi_i} f + \sum_{i\in I} \hat{\pi}(i)\hat{f}(i) \left[\log \hat{f}(i) - \log \E_{\hat{\pi}}\hat{f}\right] \\
     & = \sum_{i\in I} \hat{\pi}(i)\Ent_{\pi_i} f + \Ent_{\hat{\pi}}\hat{f}.
\end{align*}
The Dirichlet forms are
\begin{align*}
    \mathcal{L}_\pi f = & \frac{1}{2} \sum_{i,j\in I}\sum_{x,y\in\Omega} a_i(x) a_j(y)\pi(x)Q(x,y)\Psi(f(x),f(y)) \\
    = & \frac{1}{2} \sum_{i\in I}\sum_{x,y\in\Omega} a_i(x)\pi(x) a_i(y)Q(x,y)\Psi(f(x),f(y)) +   \\
    & \frac{1}{2}\sum_{i\in I}\sum_{j\neq i}\sum_{x,y\in\Omega}  a_i(x) a_j(y)\pi(x)Q(x,y)\Psi(f(x),f(y))
\end{align*}
The first term is 
\begin{align*}
    \frac{1}{2} \sum_{i\in I} \hat{\pi}(i) \sum_{x,y\in\Lambda_i} \pi_i(x) Q_i(x,y)\Psi(f(x),f(y))  = \sum_{i\in I} \hat{\pi}(i) \mathcal{L}_{\pi_i} f.
\end{align*}
 By Jensen's inequality and the convexity of $\Psi$ in all three cases,
\[\sum_{x\in\Lambda_i}\sum_{y\in\Lambda_j} \kappa_{ij}(x,y)\Psi(f(x),f(y))\geq \Psi(\hat{f}(i),\hat{f}(j)).\]
For any $(x,y,i,j)\in\Lambda_i\times\Lambda_j\times I^2$ with $x\neq y$ and $\hat{Q}(i,j)>0$, by the definition of $\chi$ we have
\[a_i(x)a_j(y)\pi(x)Q(x,y)\Psi(f(x),f(y)) \geq \chi \hat{\pi}(i)\hat{Q}(i,j) \kappa_{ij}(x,y)\Psi(f(x),f(y)).\]
This holds true even for $x=y$ or for $(i,j)$ such that $\hat{Q}(i,j)\leq 0$, so
\begin{align*}
& \sum_{x,y\in\Omega} a_i(x)a_j(y)\pi(x)Q(x,y)\Psi(f(x),f(y)) \geq \chi \hat{\pi}(i)\hat{Q}(i,j) \Psi(\hat{f}(i),\hat{f}(j))\\
\implies & \frac{1}{2}\sum_{i\in I}\sum_{j\neq i}\sum_{x,y\in\Omega}  a_i(x) a_j(y)\pi(x)Q(x,y) \Psi(f(x),f(y)) 
\end{align*}
and thus
\[\frac{1}{2}\sum_{i\in I}\sum_{j\neq i}\sum_{x,y\in\Omega}  a_i(x) a_j(y)\pi(x)Q(x,y) \Psi(f(x),f(y))
\geq \frac{1}{2}\sum_{i,j\in I}\chi \hat{\pi}(i)\hat{Q}(i,j) \Psi(\hat{f}(i),\hat{f}(j)) = \chi \mathcal{L}_{\hat{\pi}} \hat{f}\]
In summary,
\begin{align*}
    & \mathcal{R}_\pi f = \sum_{i\in I}\hat{\pi}(i)\mathcal{R}_{\pi_i}f + \mathcal{R}_{\hat{\pi}} \hat{f} \\
    & \mathcal{L}_{\pi} f \geq \sum_{i\in I} \hat{\pi}(i)\mathcal{L}_{\pi_i} f+ \chi \mathcal{L}_{\hat{\pi}} \hat{f}.
\end{align*} 
For all $\hat{\tau},\tau_i \geq 0$ such that $\mathcal{L}_{\hat{\pi}}(\hat{f})\geq \hat{\tau}\mathcal{R}_{\hat{\pi}}(\hat{f})$ and $\mathcal{L}_{\pi_i}(f)\geq \tau_i\mathcal{R}_{\pi_i}(f)$ for each $i\in I$, we have 
\begin{align*}
    \mathcal{L}_\pi(f) & \geq \sum_{i\in I}\hat{\pi}(i) \tau_i\mathcal{R}_{\pi_i}(f)
    + \chi \hat{\tau}\mathcal{R}_{\hat{\pi}}(\hat{f})  \geq \min\{\chi\hat{\tau},\min_{i\in I} \tau_i\} \mathcal{R}_{\pi}(f),
\end{align*}
and the conclusion readily follows.
\hfill\qed

\bibliography{files} 
\bibliographystyle{ieeetr}
\end{document}